\documentclass{amsart}
\usepackage{amscd}
\usepackage{amssymb}
\usepackage[T1]{fontenc}
\usepackage[utf8]{inputenc}
\usepackage{hyperref}
\usepackage[arrow,curve,matrix]{xy}
\usepackage{times}
\usepackage{graphicx}
\usepackage{color}
\usepackage{comment}

\DeclareFontFamily{OMS}{rsfs}{\skewchar\font'60}
\DeclareFontShape{OMS}{rsfs}{m}{n}{<-5>rsfs5 <5-7>rsfs7 <7->rsfs10 }{}
\DeclareSymbolFont{rsfs}{OMS}{rsfs}{m}{n}
\DeclareSymbolFontAlphabet{\scr}{rsfs}

\newcommand{\sA}{\scr{A}}
\newcommand{\sB}{\scr{B}}

\newcommand{\sF}{\scr{F}}
\newcommand{\sG}{\scr{G}}

\newcommand{\sL}{\scr{L}}

\newcommand{\sN}{\scr{N}}
\newcommand{\sO}{\scr{O}}

\newcommand{\ka}{\kappa}

\DeclareMathOperator{\codim}{codim}

\DeclareMathOperator{\Hom}{Hom}
\DeclareMathOperator{\Hilb}{Hilb}
\DeclareMathOperator{\Id}{Id}
\DeclareMathOperator{\Image}{Image}
\DeclareMathOperator{\Isom}{Isom}

\DeclareMathOperator{\Pic}{Pic}

\DeclareMathOperator{\Sym}{Sym}

\DeclareMathOperator{\Var}{Var}

\newcommand{\wtilde}{\widetilde}
\newcommand{\what}{\widehat}

\def\coh#1.#2.#3.{H^{#1}(#2,#3)}
\def\cohd#1.#2.#3.{H_{dR}^{#1}(#2,#3)}
\def\dimcoh#1.#2.#3.{h^{#1}(#2,#3)}
\def\loccoh#1.#2.#3.#4.{H^{#1}_{#2}(#3,#4)}
\def\loccohmod#1.#2.#3.{H^{#1}_{#2}(#3)}
\def\dimloccoh#1.#2.#3.#4.{h^{#1}_{#2}(#3,#4)}
\def\ses#1.#2.#3.{0  \longrightarrow  #1   \longrightarrow 
 #2 \longrightarrow #3 \longrightarrow 0} 
\def\sesshort#1.#2.#3.{0
 \rightarrow #1 \rightarrow #2 \rightarrow #3 \rightarrow 0}

\newcounter{thisthm}
\newcounter{saveenumi}

\newcommand{\iref}[1]{(\thesection.\the\value{thisthm}.\the\value{#1})}

\theoremstyle{plain}    
\newtheorem{thm}{Theorem}[section]
\newtheorem{defn}[thm]{Definition}

\numberwithin{equation}{thm}
\numberwithin{figure}{section}
\theoremstyle{plain}    
\newtheorem{cor}[thm]{Corollary}

\newtheorem{conjecture}[thm]{Conjecture}

\newtheorem{fact}[thm]{Fact}

\theoremstyle{plain}    
\newtheorem{prop}[thm]{Proposition}
\newtheorem{proclaim-special}[thm]{\specialthmname}

\theoremstyle{remark}
\newtheorem{rem}[thm]{Remark}
\newtheorem{con}[thm]{Consequence}

\newtheorem{subrem}[equation]{Remark}
\newtheorem{subclaim}[equation]{Claim} 
\newtheorem*{claim*}{Claim} 
\newtheorem{notation}[thm]{Notation}

\newtheorem{construction-defn}[thm]{Construction and Definition}

\newtheoremstyle{bozont-remark}{3pt}{3pt}%
     {}
     {}
     {\it}
     {.}
     {.5em}
     {\thmname{#1}\thmnumber{ #2}: \thmnote{\sc #3}}
\theoremstyle{bozont-remark}

\def\factor#1.#2.{\left. \raise 2pt\hbox{$#1$} \right/\hskip -2pt\raise
  -2pt\hbox{$#2$}} 
\newlength{\swidth}
\newenvironment{enumerate-p}{
  \begin{enumerate}}
  {\setcounter{equation}{\value{enumi}}\end{enumerate}}

\date{\today}

\author{Kelly Jabbusch}
\author{Stefan Kebekus}

\thanks{Kelly Jabbusch and Stefan Kebekus were supported by the
  DFG-Forschergruppe ``Classification of Algebraic Surfaces and Compact
  Complex Manifolds'' in full and in part, respectively. The work on this
  paper was finished while the authors visited the 2009 Special Year in
  Algebraic Geometry at the Mathematical Sciences Research Institute,
  Berkeley. Both authors would like to thank the MSRI for support.}
  
  \address{Kelly Jabbusch, Mathematisches Institut, Albert-Ludwigs-Universit\"at
  Freiburg, Eckerstra{\ss}e 1, 79104 Freiburg, Germany}
\email{\href{mailto:kelly.jabbusch@math.uni-freiburg.de}{kelly.jabbusch@math.uni-freiburg.de}}
  
\address{Stefan Kebekus, Mathematisches Institut, Albert-Ludwigs-Universit\"at
  Freiburg, Eckerstra{\ss}e 1, 79104 Freiburg, Germany}
\email{\href{mailto:stefan.kebekus@math.uni-freiburg.de}{stefan.kebekus@math.uni-freiburg.de}}
\urladdr{\href{http://home.mathematik.uni-freiburg.de/kebekus}{http://home.mathematik.uni-freiburg.de/kebekus}}

\title{Positive sheaves of differentials coming from coarse moduli spaces}

\date{\today}
\begin{document}

\begin{abstract}
  Consider a smooth projective family of canonically polarized complex
  manifolds over a smooth quasi-projective complex base $Y^\circ$, and suppose
  the family is non-isotrivial.  If $Y$ is a smooth compactification of
  $Y^\circ$, such that $D:=Y \setminus Y^\circ$ is a simple normal crossing
  divisor, then we can consider the sheaf of differentials with logarithmic
  poles along $D$.  Viehweg and Zuo have shown that for some $m > 0$, the
  $m^{\rm th}$ symmetric power of this sheaf admits many sections. More
  precisely, the $m^{\rm th}$ symmetric power contains an invertible sheaf
  whose Kodaira-Iitaka dimension is at least the variation of the family.  We
  refine this result and show that this ``Viehweg-Zuo sheaf'' comes from the
   coarse moduli space associated to the given family, at least generically.

   As an immediate corollary, if $Y^\circ$ is a surface, we see that the
   non-isotriviality assumption implies that $Y^\circ$ cannot be special in the
   sense of Campana.
 \end{abstract}

 \maketitle
 \tableofcontents

\section{Introduction and statement of main result}
\label{sec:VZ}

\subsection{Introduction}

Throughout this paper, we consider a smooth projective family $f^\circ:
X^\circ \to Y^\circ$ of canonically polarized complex manifolds, of relative
dimension $n$, over a smooth complex quasi-projective base. We assume that the
family is not isotrivial, and let $\mu: Y^\circ \to \mathfrak{M}$ be the
associated map to the coarse moduli space, whose existence is shown, e.g.~in
\cite[Thm.~1.11]{V95}. We fix a smooth projective compactification $Y$ of
$Y^\circ$ such that $D := Y \setminus Y^\circ$ is a divisor with simple normal
crossings.

In this setup, Viehweg and Zuo have shown the following fundamental result
concerning the existence of logarithmic pluri-differentials on $Y^\circ$.

\begin{thm}[Existence of pluri-differentials on $Y$, \protect{\cite[Thm.~1.4(i)]{VZ02}}]\label{thm:VZ}
  There exists a number $m > 0$ and an invertible sheaf $\sA \subseteq \Sym^m
  \Omega^1_Y (\log D)$ whose Kodaira-Iitaka dimension is at least the
  variation of the family, $\kappa(\sA) \geq \Var(f^\circ)$. \qed
\end{thm}

The ``Viehweg-Zuo'' sheaf $\sA$ was crucial in the study of hyperbolicity
properties of manifolds that appear as bases of families of maximal variation
and has been used to show that any minimal model program of the pair $(Y, D)$
factors the moduli map, \cite{KK08, KK08b, KK08c}, see also the survey
\cite{KS06}. In spite of its importance, little is known about further
properties of the sheaf $\sA$.  For example, it is unclear to us how the
Viehweg-Zuo construction behaves under base change. The goal of this short
note is to refine the result of Viehweg and Zuo somewhat, and show that the
Viehweg-Zuo sheaf $\sA$ comes from the coarse moduli space $\mathfrak{M}$, at
least generically. A precise statement is given in
Theorem~\ref{thm:VZimproved} below.

Theorem~\ref{thm:VZimproved} directly relates to a conjecture of Campana. In
\cite[Conj.~12.19]{Cam07} Campana conjectured that the assumption ``$f^\circ$
not isotrivial'' immediately implies that the base manifold $Y^\circ$ is not
\emph{special}. In other words, any family of canonically polarized varieties
over a special base manifold is necessarily isotrivial. In the case where
$Y^\circ$ is a surface, the conjecture is claimed in
\cite[Thm.~12.20]{Cam07}. However, we had difficulties following the proof.
We will show in Section~\ref{sec:special} that Campana's conjecture in
dimension two is an immediate corollary to
Theorem~\ref{thm:VZimproved}. Using a more advanced line of argumentation,
Campana's conjecture in dimension three can also be deduced. Details will
appear in a forthcoming paper.
 
Throughout the present paper we work over the field of complex numbers.

\subsection{Statement of the main result}

Roughly speaking, the main result of this paper is that the
Viehweg-Zuo sheaf comes from the coarse moduli space $\mathfrak M$. To
formulate the statement precisely, we use the following notation.

\begin{notation}\label{not:introB}
  Consider the subsheaf $\sB \subseteq \Omega^1_Y(\log D)$, defined on
  presheaf level as follows: if $U \subset Y$ is any open set and $\sigma \in
  \Gamma\bigl(U,\, \Omega^1_Y (\log D) \bigr)$ any section, then $\sigma \in
  \Gamma\bigl(U,\, \sB \bigr)$ if and only if the restriction
  $\sigma|_{U'}$ is in the image of the differential map
  $$
  d\mu|_{U'} : \mu^* \bigl( \Omega^1_{\mathfrak M}\bigr)|_{U'} \longrightarrow \Omega^1_{U'},
  $$
  where $U' \subseteq U\cap Y^\circ$ is the open subset where the moduli map
  $\mu$ has maximal rank.
\end{notation}
\begin{rem}
  By construction, it is clear that the sheaf $\sB$ is a saturated subsheaf of
  $\Omega^1_Y (\log D)$. We say that $\sB$ is the saturation of $\Image(d\mu)$
  in $\Omega^1_Y(\log D)$.
\end{rem}

With this notation, the main result of the paper is then formulated
as follows.

\begin{thm}[Refinement of the Viehweg-Zuo Theorem~\ref{thm:VZ}]\label{thm:VZimproved}
  There exists a number $m > 0$ and an invertible subsheaf $\sA \subseteq
  \Sym^m \sB$ whose Kodaira-Iitaka dimension is at least the variation of the
  family, $\kappa(\sA) \geq \Var(f^\circ)$.
\end{thm}

\subsection{Outline of the paper}

We begin the proof of Theorem~\ref{thm:VZimproved} in
Section~\ref{sec:synopsis} with a summary of Viehweg-Zuo's proof of
Theorem~\ref{thm:VZ}. Using the notation and results of
Section~\ref{sec:synopsis}, a proof of Theorem~\ref{thm:VZimproved} is given
in Section~\ref{sec:pfVZimproved}. We end this paper with
Section~\ref{sec:special}, where we briefly recall Campana's notion of a
special logarithmic pair, give the precise statement of his conjecture and
give an extremely short proof for families over surfaces.

\subsection*{Acknowledgments}

Campana's Conjecture~\ref{conj:campana} was brought to our attention during
the 2007 Levico conference in Algebraic Geometry. We would like to thank
Fr\'ed\'eric Campana and Eckart Viehweg for a number of discussions on the
subject.

\section{Review of Viehweg-Zuo's proof of Theorem~\ref*{thm:VZ}}
\label{sec:synopsis}

To prepare for the proof of Theorem~\ref{thm:VZimproved}, we give a very brief
synopsis of Viehweg-Zuo's proof of Theorem~\ref{thm:VZ}, covering only the
material used in the proof of Theorem~\ref{thm:VZimproved}.  The reader who is
interested in a detailed understanding is referred to the original paper
\cite{VZ02} and to the survey \cite{Viehweg01}. The overview contained in this
section and the facts outlined in Section~\ref{ssec:facts} can perhaps serve
as a guideline to the original references.

\subsection{Setup of notation}
\label{ssec:8A}

Throughout the present Section~\ref{sec:synopsis}, we choose a smooth
projective compactification $X$ of $X^\circ$ such that the following holds:
\begin{enumerate-p}
\item The difference $\Delta := X \setminus X^\circ$ is a divisor with simple
  normal crossings.
\item The morphism $f^\circ$ extends to a projective morphism $f: X \to Y$.
  \setcounter{saveenumi}{\theenumi}
\end{enumerate-p}
It is then clear that $\Delta = f^{-1}(D)$ set-theoretically. Removing a
suitable subset $S \subset Y$ of codimension $\codim_Y S \geq 2$, the
following will then hold automatically on $Y' := Y\setminus S$ and $X' := X
\setminus f^{-1}(S)$, respectively.
\begin{enumerate-p}
\setcounter{enumi}{\thesaveenumi}
\item The restricted morphism $f' := f|_{X'}$ is flat.
\item The divisor $D' := D\cap Y'$ is smooth.
\item The divisor $\Delta' := \Delta \cap X'$ is a relative normal crossing
  divisor, i.e. a normal crossing divisor whose components and all their
  intersections are smooth over the components of $D'$.
\end{enumerate-p}
In the language of Viehweg-Zuo, \cite[Def~2.1(c)]{VZ02}, the restricted
morphism $f' : X' \to Y'$ is a ``good partial compactification of $f^\circ$''.

\begin{rem}[Restriction to a partial compactification]\label{rem:extoS}
  Let $\sG$ be a locally free sheaf on $Y$, and let $\sF' \subseteq \sG|_{Y'}$
  be an invertible subsheaf. Since $\codim_Y S \geq 2$, there exists a unique
  extension of the sheaf $\sF'$ to an invertible subsheaf $\sF \subseteq \sG$
  on $Y$. Furthermore, the restriction map $\Gamma \bigl( Y,\, \sF \bigr) \to
  \Gamma \bigl( Y',\, \sF' \bigr)$ is an isomorphism. In particular, the
  notion of Kodaira-Iitaka dimension makes sense for the sheaf $\sF'$, and
  $\kappa(\sF') = \kappa(\sF)$.
\end{rem}

\subsection{Construction of the $\boldsymbol{\tau^0_{p,q}}$}

The starting point of the Viehweg-Zuo construction is the standard sequence of
logarithmic differentials associated to the flat morphism $f'$,
\begin{equation}\label{eq:lds}
  0 \to (f')^*\Omega^1_{Y'}(\log D') \to \Omega^1_{X'}(\log \Delta') \to
  \Omega^1_{X'/Y'}(\log \Delta') \to 0,
\end{equation}
where $\Omega^1_{X'/Y'}(\log \Delta')$ is locally free. It is a standard fact,
\cite[II, Ex.~5.16]{Ha77}, that Sequence~\eqref{eq:lds} defines a filtration
of the $p^{\rm th}$ exterior power,
$$
\Omega^p_{X'}(\log \Delta') = F^0 \supseteq F^1 \supseteq \cdots \supseteq F^p
\supseteq F^{p+1} = \{0\},
$$
with $F^r/F^{r+1} \cong (f')^*\bigl( \Omega^r_{Y'}(\log D') \bigr) \otimes
\Omega^{p-r}_{X'/Y'}(\log \Delta')$. Taking the sequence
$$
0 \longrightarrow F^1 \longrightarrow F^0 \longrightarrow \factor F^0.F^1.
\longrightarrow 0
$$
modulo $F^2$, we obtain
\begin{equation}\label{eq:amt}
  0 \to (f')^*\bigl( \Omega^1_{Y'}(\log D')\bigr) \otimes \Omega^{p-1}_{X'/Y'}(\log \Delta')
  \to \factor F^0.F^2. \to \Omega^p_{X'/Y'}(\log \Delta') \to 0.
\end{equation}
Setting $\sL := \Omega^n_{X'/Y'} (\log \Delta')$, twisting
Sequence~\eqref{eq:amt} with $\sL^{-1}$ and pushing down, the connecting
morphisms of the associated long exact sequence give maps
$$
\tau^0_{p,q}: F^{p,q} \longrightarrow F^{p-1,q+1} \otimes \Omega_{Y'}^1(\log
D'),
$$
where $F^{p,q}:= \factor R^q f'_*( \Omega_{X'/Y'}^p (\log \Delta') \otimes
\sL^{-1}).\text{torsion}.$. Set $\sN_0^{p,q} := \ker (\tau^0_{p,q})$.

\subsection{Alignment of the $\boldsymbol{\tau^0_{p,q}}$}

The morphisms $\tau^0_{p,q}$ and $\tau^0_{p-1,q+1}$ can be composed if we
tensor the latter with the identity morphism on $\Omega_{Y'}^1(\log D')$. More
specifically, we consider the following morphisms,
$$
\underbrace{\tau^0_{p,q} \otimes \Id_{\Omega^1_{Y'}(\log D')^{\otimes q}}}_{=:
  \tau_{p,q}} : F^{p,q} \otimes \bigl(\Omega^1_{Y'}(\log D')\bigr)^{\otimes q}
\to F^{p-1,q+1} \otimes \bigl( \Omega_{Y'}^1(\log D') \bigr)^{\otimes q+1},
$$
and their compositions
\begin{equation}\label{eq:tk}
  \underbrace{\tau_{n-k+1,k-1} \circ \cdots \circ \tau_{n-1,1} \circ \tau_{n,0}}_{=:
    \tau^k} : F^{n,0} \to F^{n-k, k} \otimes \bigl(\Omega_{Y'}^1(\log D')
  \bigr)^{\otimes k}.
\end{equation}

\subsection{Fundamental facts about $\boldsymbol{\tau^k}$ and $\boldsymbol{\sN_0^{p,q}}$}
\label{ssec:facts}

Theorem~\ref{thm:VZ} is shown by relating the morphism $\tau^0_{p,q}$ with the
structure morphism of a Higgs-bundle coming from the variation of Hodge
structures associated with the family $f^\circ$. Viehweg's positivity results
of push-forward sheaves of relative differentials, as well as Zuo's results on
the curvature of kernels of generalized Kodaira-Spencer maps are the main
input here. Rather than recalling the complicated line of argumentation, we
simply state two central results from the argumentation of \cite{VZ02}.

\begin{fact}[Factorization via symmetric differentials, \protect{\cite[Lem.~4.6]{VZ02}}]\label{fact:VZ46}
  For any $k$, the morphism $\tau^k$ factors via the symmetric differentials
  $\Sym^k \Omega_{Y'}^1(\log D') \subseteq \bigl(\Omega_{Y'}^1(\log D')
  \bigr)^{\otimes k}$. More precisely, the morphism $\tau^k$ takes its image
  in $F^{n-k, k} \otimes \Sym^k \Omega_{Y'}^1(\log D')$.  \qed
\end{fact}

\begin{con}
  Using Fact~\ref{fact:VZ46} and the observation that $F^{n,0} \cong
  \sO_{Y'}$, we can therefore view $\tau^k$ as a morphism
  $$
  \tau^k : \sO_{Y'} \longrightarrow F^{n-k, k} \otimes \Sym^k
  \Omega_{Y'}^1(\log D').
  $$
\end{con}

While the proof of Fact~\ref{fact:VZ46} is rather elementary, the following
deep result is at the core of Viehweg-Zuo's argument.

\begin{fact}[Negativity of $\sN_0^{p,q}$, \protect{\cite[Claim~4.8]{VZ02}}]\label{fact:VZ48}
  Given any numbers $p$ and $q$, there exists a number $k$ and an invertible
  sheaf $\sA' \in \Pic(Y')$ of Kodaira-Iitaka dimension $\kappa(\sA') \geq
  \Var(f^0)$ such that $(\sA')^* \otimes \Sym^k \bigl( (\sN_0^{p,q})^* \bigr)$
  is generically generated. \qed
\end{fact}

\subsection{End of proof}
\label{ssec:eopVZ}

To end the sketch of proof, we follow \cite[p.~315]{VZ02} almost verbatim. By
Fact~\ref{fact:VZ48}, the trivial sheaf $F^{n,0} \cong \sO_{Y'}$ cannot lie in
the kernel $\sN_0^{n,0}$ of $\tau^1 = \tau^0_{n,0}$. We can therefore set $1
\leq m$ to be the largest number  with $\tau^{m}(F^{n,0}) \not = \{0\}$.
Since $m$ is maximal, $\tau^{m+1} = \tau_{n-m,m} \circ \tau^{m} \equiv 0$ and
$$
\Image(\tau^{m}) \subseteq \ker(\tau_{n-m,m}) = \sN_0^{n-m,m}\otimes \Sym^m
\Omega^1_{Y'}(\log D').
$$
In other words, $\tau^{m}$ gives a non-trivial map
$$
\tau^{m}: \sO_{Y'} \cong F^{n,0} \longrightarrow \sN_0^{n-m,m}\otimes \Sym^m
\Omega^1_{Y'}(\log D').
$$
Equivalently, we can view $\tau^{m}$ as a non-trivial map
\begin{equation}\label{eq:end}
  \tau^{m}: (\sN_0^{n-m,m})^*  \longrightarrow \Sym^m \Omega^1_{Y'}(\log D').
\end{equation}
By Fact~\ref{fact:VZ48}, there are many morphisms $\sA' \to \Sym^k \bigl(
(\sN_0^{n-m, m})^* \bigr)$, for $k$ large enough. Together
with~\eqref{eq:end}, this gives a non-zero morphism $\sA' \to \Sym^{k\cdot m}
\Omega^1_{Y'}(\log D')$.

We have seen in Remark~\ref{rem:extoS} that the sheaf $\sA' \subseteq
\Sym^{k\cdot m} \Omega^1_{Y'}(\log D')$ extends to a sheaf $\sA \subseteq
\Sym^{k\cdot m} \Omega^1_Y(\log D)$ with $\kappa(\sA) = \kappa(\sA') \geq
\Var(f^\circ)$. This ends the proof of Theorem~\ref{thm:VZ}. \qed

\section{Proof of Theorem~\ref*{thm:VZimproved}}
\label{sec:pfVZimproved}

\subsection{Setup and assumptions}

The proof of Theorem~\ref{thm:VZimproved} makes use of essentially all results
explained in Section~\ref{sec:synopsis}. Since the assumptions of
Theorems~\ref{thm:VZ} and \ref{thm:VZimproved} agree, we maintain the full
setup and all notation introduced in Section~\ref{sec:synopsis}.

\subsection{Reduction to a study of the $\boldsymbol{\tau^0_{p,q}}$}

The construction outlined in Section~\ref{sec:synopsis} essentially says that
the sheaf $\sA$ constructed by Viehweg-Zuo is a symmetric product of the image
sheaves of the $\tau^0_{p,q}$. The precise statement is the following.

\begin{prop}\label{prop:VZConstruction}
  To prove Theorem~\ref{thm:VZimproved}, it suffices to show that
  \begin{equation}\label{eq:tau}
    \Image(\tau^0_{p,q}) \subseteq F^{p-1,q+1} \otimes \sB'
  \end{equation}
  for all $p$ and $q$, where $\sB' := \sB|_{Y'}$ and $\sB \subseteq
  \Omega^1_Y(\log D)$ is the sheaf defined in Notation~\ref{not:introB}.
\end{prop}

\begin{subrem}\label{rem:open}
  Since $\sB$ is saturated, it is enough to check inclusion \eqref{eq:tau} on
  an open set.
\end{subrem}

\begin{proof}[Proof of Proposition~\ref{prop:VZConstruction} ]
  If~\eqref{eq:tau} holds, the image of the morphisms $\tau^k$ defined
  in~\eqref{eq:tk} lies in $F^{n-k,k} \otimes \bigl( \sB' \bigr)^{\otimes
    k}$. Furthermore, by Fact~\ref{fact:VZ46},
  $$
  \Image(\tau^k) \subseteq F^{n-k,k} \otimes \Sym^k \sB'.
  $$
  If we chose the number $m$ as in Section~\ref{ssec:eopVZ} above, the image
  of $\tau^{m}$ is then contained in $\sN^{n-m,m}_0 \otimes \Sym^m \sB'$,
  and $\tau^{m}$ can be seen as a non-trivial map
  $$
  \tau^{m} : (\sN_0^{n-m,m})^*  \longrightarrow \Sym^m \sB'.
  $$
  As in Section~\ref{ssec:eopVZ}, we obtain a map $\sA' \to \Sym^{k\cdot m}
  \sB'$, with $\kappa(\sA') \geq \Var(f^\circ)$, and Remark~\ref{rem:extoS}
  gives the extension to a sheaf $\sA \subset \Sym^{k\cdot m} \sB$, with
  $\kappa(\sA) = \kappa(\sA')$.
\end{proof}

\subsection{Proof of Inclusion~\eqref{eq:tau} in a simple case}

It remains to check Inclusion~\eqref{eq:tau}. Before tackling the problem in
general, we consider a trivial case first.

\begin{prop}\label{prop:incformaxvar}
  If the variation of $f^\circ$ is maximal, i.e.~$\Var(f^\circ) = \dim
  Y^\circ$, then Inclusion~\eqref{eq:tau} holds.
\end{prop}
\begin{proof}
  If the variation of $f^\circ$ is maximal, then the moduli map $Y^\circ \to
  \mathfrak M$ is generically finite onto the closure of its image. In
  particular, the sheaf $\sB$ introduced in Notation~\ref{not:introB} equals
  $\Omega^1_Y(\log D)$.  Inclusion~\eqref{eq:tau} is therefore trivially
  satisfied.
\end{proof}

\subsection{Comparing families with respect to Inclusion~\eqref{eq:tau}}

Given two families, one the pull-back of the other via a dominant morphism, an
elementary comparison of the morphisms $\tau^0_{p,q}$ associated with the
families shows that Inclusion~\eqref{eq:tau} holds for one of the families if
and only if it also holds for the other. We will later use the following
Comparison Proposition to show that the Viehweg-Zuo sheaf of a family
essentially only depends on the image of the base in the coarse moduli space,
and not so much on the family itself.

\begin{prop}[Comparison Proposition]\label{prop:comparison}
  Consider a Cartesian diagram of smooth projective families of $n$-dimensional
  canonically polarized manifolds over smooth quasi-projective base manifolds,
  as follows
  $$
  \xymatrix{
    \hat X^\circ \ar[rr]^{\Gamma} \ar[d]_{\hat f^\circ} &&
    \tilde X^\circ \ar[d]^{\tilde f^\circ} \\
    \hat Y^\circ \ar[rr]^{\gamma}_{\text{dominant}} && \tilde
    Y^\circ.  }
  $$
  Let $\hat f' : \hat X' \to \hat Y'$ and $\tilde f' : \tilde X' \to \tilde
  Y'$ be two good partial compactifications, in the sense introduced in
  Section~\ref{ssec:8A}. Then Inclusion~\eqref{eq:tau} holds for $\hat f'$ if
  and only if it holds for $\tilde f'$.
\end{prop}
\begin{proof}
  We have noted in Remark~\ref{rem:open} that it suffices to check
  Inclusion~\eqref{eq:tau} on an open subset. In particular, it suffices to
  consider the restrictions of the morphisms $\tau^0_{p,q}$ and of all
  relevant sheaves to $\hat Y^\circ$ and $\tilde Y^\circ$. This greatly
  simplifies notation because the logarithmic boundary terms do not appear in
  the restrictions, and we can write, e.g., $\Omega^p_{\tilde Y^\circ}$
  instead of the more complicated $\Omega^p_{\tilde Y'}(\log \tilde D')$.

  Shrinking $\hat Y^\circ$ and $\tilde Y^\circ$ further, if necessary, we may
  assume without loss of generality that $\gamma$ is surjective and smooth. We
  may also assume that the moduli map $\tilde \mu : \tilde Y^\circ \to
  \mathfrak M$ has maximal rank. By assumption, the moduli map $\hat \mu :
  \hat Y^\circ \to \mathfrak M$ is the composition $\hat \mu = \tilde \mu
  \circ \gamma$.

  As in Section~\ref{sec:synopsis}, we need to discuss the connecting
  morphisms $\tau^0_{p,q}$ on $\hat Y^\circ$ and on $\tilde Y^\circ$,
  respectively. For clarity of notation we indicate the relevant space by
  indexing all morphisms and sheaves with either a hat or a tilde.  That way,
  we write
  $$
  \hat \tau^0_{p,q} : \hat F^{p,q} \to \hat F^{p-1,q+1} \otimes \Omega^1_{\hat
    Y^\circ} \text{\,\, and \,\,} \tilde \tau^0_{p,q} : \tilde F^{p,q} \to
  \tilde F^{p-1,q+1}\otimes \Omega^1_{\tilde Y^\circ},
  $$
  where $\hat F^{p,q} := R^q\hat f^\circ_*\bigl( \Omega^p_{\hat X^\circ/\hat
    Y^\circ} \otimes (\Omega^n_{\hat X^\circ/\hat Y^\circ})^{-1} \bigr)$ and
  the sheaf $\tilde F^{p,q}$ on $\tilde Y^\circ$ is defined
  analogously. Finally, set
  $$
\what \sB^\circ := \Image \bigl( d\hat \mu : \hat
  \mu^*(\Omega^1_{\mathfrak M}) \to \Omega^1_{\hat Y^\circ}\bigr)
  \text{\,\, and \,\,}
  \wtilde \sB^\circ := \Image\bigl(d\tilde \mu : \tilde
  \mu^*(\Omega^1_{\mathfrak M}) \to \Omega^1_{\tilde Y^\circ} \bigr). 
  $$
  Since $\gamma$ is smooth and the moduli map $\tilde \mu$ has maximal rank,
  $\hat \mu$ also has maximal rank, and both $\wtilde \sB^\circ$ and $\what
  \sB^\circ$ are saturated in $\Omega^1_{\tilde Y^\circ}$ and $\Omega^1_{\hat
    Y^\circ}$, respectively. Better still, the differential $d\gamma :
  \gamma^*(\Omega^1_{\tilde Y^\circ}) \to \Omega^1_{\hat Y^\circ}$ induces an
  isomorphism
  \begin{equation}\label{eq:dgB}
    d\gamma : \gamma^* (\wtilde \sB^\circ) \xrightarrow{\,\, \cong \,\,}
    \what \sB^\circ.
  \end{equation}
  Since $\wtilde \sB^\circ$ and $\what \sB^\circ$ are saturated, to prove
  Proposition~\ref{prop:comparison} it is equivalent to show that
  \begin{equation}\label{eq:Bequiv}
    \Image(\hat \tau^0_{p,q}) \subseteq \hat F^{p-1,q+1} \otimes \what
    \sB^\circ \Longleftrightarrow \Image(\tilde \tau^0_{p,q}) \subseteq
    \tilde F^{p-1,q+1} \otimes \wtilde \sB^\circ.
  \end{equation}
  To prove~\eqref{eq:Bequiv}, we aim to identify the sheaves $\hat F^{p,q}$
  and $\gamma^* \bigl(\tilde F^{p,q}\bigr)$ and show that images of the $\hat
  \tau^0_{p,q}$ are naturally identified with the pull-backs of the images of
  $\tilde \tau^0_{p,q}$. For a precise statement, recall that there are
  isomorphisms $\Gamma^* \bigl( \Omega^p_{\tilde X^\circ/\tilde Y^\circ}
  \bigr) \cong \Omega^p_{\hat X^\circ/\hat Y^\circ}$ for all $p$. Since taking
  cohomology commutes with flat base change, \cite[III Prop.~9.3]{Ha77}, we
  obtain isomorphisms
  $$
  \iota^{p,q} : \gamma^*\bigl(\tilde F^{p,q}\bigr) \xrightarrow{\,\, \cong
    \,\,} \hat F^{p,q}
  $$
  for all $p$ and $q$. Tensoring $\iota^{p,q}$ with the differential $d\gamma
  : \gamma^* \bigl( \Omega^1_{\tilde Y^\circ} \bigr) \to \Omega^1_{\hat
    Y^\circ}$ gives a map
  \begin{equation}\label{eq:smnqi}
    \iota^{p,q} \otimes d\gamma : \gamma^* \bigl( \tilde F^{p,q} \otimes
    \Omega^1_{\tilde Y^\circ} \bigr) \longrightarrow \hat F^{p,q}
    \otimes \Omega^1_{\hat Y^\circ}.
  \end{equation}
  Equivalence~\eqref{eq:Bequiv}, and hence Proposition~\ref{prop:comparison},
  is an immediate consequence of the Isomorphism~\eqref{eq:dgB} and of the
  following claim.

  \begin{subclaim}\label{sclaim:putaupq}
    Given any numbers $p$ and $q$, the sheaf morphism~\eqref{eq:smnqi} induces
    an isomorphism between the image of $\hat \tau^0_{p,q}$ and the pull-back of
    the image of $\tilde \tau^0_{p,q}$,
    $$
    \iota^{p,q} \otimes d\gamma : \gamma^* \bigl( \Image(\tilde \tau^0_{p,q})
    \bigr) \xrightarrow{\,\, \cong \,\,} \Image(\hat \tau^0_{p,q}).
    $$
  \end{subclaim}

  It remains to prove Claim~\ref{sclaim:putaupq}.  Observe that
  Claim~\ref{sclaim:putaupq} follows trivially from the definitions of $\tilde
  \tau^0_{p,q}$ and $\hat \tau^0_{p,q}$ if we are in the simple case where
  $\hat Y^0$ is a product, say $\hat Y^\circ \cong \tilde Y^\circ \times
  \tilde Z^\circ$, and where $\gamma$ is the projection to the first
  factor. Locally in the analytic topology, however, any smooth morphism looks
  like the projection morphism of a product.  Since Claim~\ref{sclaim:putaupq}
  can be checked locally analytically, this proves the claim and ends the
  proof of Proposition~\ref{prop:comparison}.
\end{proof}

\subsection{End of proof of Theorem \ref{thm:VZimproved}}

To complete the proof of Theorem~\ref{thm:VZimproved}, we compare our original
family to one that is of maximal variation.  The starting point is the
existence of a universal family on a finite cover.

\begin{thm}[\protect{Existence of a universal family on a finite cover, \cite[Prop.~2.7]{Kollar90}, see also \cite[Thm.~9.25]{V95}}]
  Let $\mathfrak M' \subseteq \mathfrak M$ be the reduced irreducible
  component that contains the image of $Y^\circ$. Then there exists a reduced
  normal scheme $\wtilde{\mathfrak M}$, a finite and surjective morphism
  $\gamma : \wtilde{\mathfrak M} \to \mathfrak M'$ and a family of canonically
  polarized varieties $u: \tilde U \to \wtilde{\mathfrak M}$ such that
  $\gamma$ is precisely the moduli map associated with the family $u$. \qed
\end{thm}

Let $Z' \subseteq Y^\circ \times_{\mathfrak M} \wtilde{\mathfrak M}$ be an
irreducible component of the fiber product that dominates $Y^\circ$, and let
$Z$ be a desingularization of $Z'$. Setting $X^\circ_Z := X^\circ
\times_{Y^\circ} Z$ and $\tilde U_Z := \tilde U \times_{\wtilde{\mathfrak M}}
Z$, we obtain two linked Cartesian diagrams, as follows
$$
\xymatrix{ X^\circ \ar[d]_{f^\circ} && X^\circ_Z \ar[d]_{f^\circ_Z} \ar[ll] &
  \tilde U_Z \ar[d]_{u_Z}^{\txt{\tiny family of\\\tiny max.~var.}} \ar[rr] &&
  \tilde U \ar[d]_{u} \\
  Y^\circ && Z \ar[ll]_{\txt{\tiny dominant}} \ar@{=}[r]& Z
  \ar[rr]^{\txt{\tiny gen.~finite}} && \wtilde{\mathfrak
    M} \ar[rr]^{\txt{\tiny finite, surjective}}_{\txt{\tiny moduli map of $u$}} &&
  \mathfrak M.}
$$
Here $f^\circ_Z$ and $u_Z$ are two families of canonically polarized varieties
that are not necessarily isomorphic, but induce the same moduli map $Z \to
\mathfrak M$. Since for any point $z \in Z$, the fibers $(f^\circ_Z)^{-1}(z)$
and $u_Z^{-1}(z)$ are isomorphic, the scheme of $Z$-isomorphisms,
$$
I' := \Isom_Z \bigl( X^\circ_Z,\, \tilde U_Z \bigr) \subseteq \Hom_Z \bigl(
X^\circ_Z,\, \tilde U_Z \bigr) \subseteq \Hilb_Z \bigl( X^\circ_Z\times_Z
\tilde U_Z \bigr)
$$
surjects onto $Z$. Since all fibers $(f^\circ_Z)^{-1}(z) \cong u_Z^{-1}(z)$
are canonically polarized manifolds and have only finitely many automorphisms,
the natural map $I' \to Z$ is quasi-finite. Let $I$ be a desingularization of
a component of $I'$ that dominates $Z$. Recall that taking $\Hilb$, $\Hom$ and
$\Isom$ commutes with base change. In particular, we have an isomorphism of
$I$-schemes,
$$
\Isom_I \bigl( X^\circ_Z \times_Z I,\, \tilde U_Z \times_Z I \bigr) \cong
\Isom_Z \bigl( X^\circ_Z,\, \tilde U_Z \bigr) \times_Z I.
$$
Looking at the right hand side, it is clear that there exists a section $I \to
\Isom_I \bigl( X^\circ_Z \times_Z I,\, \tilde U_Z \times_Z I \bigr)$, i.e., an
isomorphism of $I$-schemes, $X^\circ_Z \times_Z I \cong \tilde U_Z \times_Z
I$. In summary, we obtain a diagram as follows,
$$
\xymatrix{
  X^\circ \ar[d]_{f^\circ} && X^\circ_Z \times_Z I \ar[d]_{f^\circ_I} \ar[ll]
  \ar@{<->}[r]^{\cong} & \tilde U_Z \times_Z I \ar[d]^{u_I} \ar[rr] &&
  \tilde U_Z \ar[d]_{u_Z}^{\txt{\tiny family of\\\tiny max.~var.}} \\
  Y^\circ && I \ar[ll]_{\gamma_Y}^{\txt{\tiny dominant}} \ar@{=}[r]& I
  \ar[rr]^{\gamma_Z}_{\txt{\tiny dominant}} && Z. }
$$
In this situation, Proposition~\ref{prop:incformaxvar} asserts that
Inclusion~\eqref{eq:tau} holds for the family $u_Z$, which is of maximal
variation. Since $\gamma_Z$ is dominant, and since it suffices to prove
Inclusion~\eqref{eq:tau} on an open subset, the Comparison
Proposition~\ref{prop:comparison} applies to show that
Inclusion~\eqref{eq:tau} holds for the family $u_I$ or equivalently for the
family $f^\circ_I$. Another application of the Comparison
Proposition~\ref{prop:comparison} to the morphism $\gamma_Y$ then shows that
Inclusion~\eqref{eq:tau} holds for the family $f^\circ$.
Theorem~\ref{thm:VZimproved} then follows from
Proposition~\ref{prop:VZConstruction}. \qed

\section{Application of Theorem~\ref*{thm:VZimproved} to families over special surfaces}
\label{sec:special}

As an immediate corollary to Theorem~\ref{thm:VZimproved}, we see that any
smooth projective family of canonically polarized manifolds over a special
surface $Y^\circ$ is isotrivial, as conjectured by Campana.  We first recall
the precise definition of a special logarithmic pair below, taking the
classical Bogomolov-Sommese vanishing theorem as our starting point.

\begin{thm}[\protect{Bogomolov-Sommese vanishing, \cite[Sect.~6]{EV92}}]\label{thm:classBSv}
  Let $Y$ be a smooth projective variety and $D \subset Y$ a reduced, possibly
  empty divisor with simple normal crossings. If $p \leq \dim Y$ is any number
  and $\sA \subseteq \Omega^p_Y(\log D)$ any invertible subsheaf, then the
  Kodaira-Iitaka dimension of $\sA$ is at most $p$, i.e., $\kappa(\sA) \leq
  p$. \qed
\end{thm}

In a nutshell, we say that a pair $(Y, D)$ is special if the inequality in the
Bogomolov-Sommese vanishing theorem is always strict.

\begin{defn}[Special logarithmic pair]\label{def:speciallog}
  In the setup of Theorem~\ref{thm:classBSv}, a pair $(Y, D)$ is called
  \emph{special} if the strict inequality $\kappa(\sA) < p$ holds for all $p$
  and all invertible sheaves $\sA \subseteq \Omega^p_Y(\log D)$. A smooth,
  quasi-projective variety $Y^\circ$ is called special if there exists a
  smooth compactification $Y$ such that $D:= Y \setminus Y^\circ$ is a divisor
  with simple normal crossings and such that the pair $(Y,D)$ is
  special.
\end{defn}

\begin{rem}
  If $Y^\circ$ is a smooth, quasi-projective variety and if $(Y_1, D_1)$ and
  $(Y_2, D_2)$ are two smooth compactifications with snc boundary divisors, as
  in Definition~\ref{def:speciallog}, then an elementary computation shows
  that the pair $(Y_1, D_1)$ is special if and only if $(Y_2, D_2)$ is
  special. Specialness can thus be checked on any snc compactification.
\end{rem}

With this notation in place, Campana has conjectured the following.

\begin{conjecture}[\protect{Generalization of Shafarevich Hyperbolicity, \cite[Conj.~12.19]{Cam07}}]\label{conj:campana}
  Let $f: X^\circ \to Y^\circ$ be a smooth family of canonically polarized
  varieties over a smooth quasi-projective base. If $Y^\circ$ is special, then
  the family $f^\circ$ is isotrivial.
\end{conjecture}

As mentioned in the Introduction, in the case where $Y^\circ$ is a surface,
Conjecture~\ref{conj:campana} is claimed in \cite[Thm.~12.20]{Cam07}. However,
we had difficulties following the proof, and offer a new proof, which is an
immediate corollary to Theorem~\ref{thm:VZimproved}.

\begin{cor}[Campana's conjecture in dimension two]
  Conjecture~\ref{conj:campana} holds if $\dim Y^\circ = 2$.
\end{cor}
\begin{proof}
  We maintain the notation of Conjecture~\ref{conj:campana} and let $f:
  X^\circ \to Y^\circ$ be a smooth family of canonically polarized varieties
  over a smooth quasi-projective base, with $Y^\circ$ a special surface.
  Since $Y^\circ$ is special, it is not of log general type, and hence by
  \cite[Thm.~1.1]{KK08c}, $\Var(f^\circ) < 2$.  Suppose $\Var(f^\circ)=1$ and
  choose a compactification $(Y, D)$ as in Definition~\ref{def:speciallog},
  then by Theorem~\ref{thm:VZimproved} there exists a number $m>0$ and an
  invertible subsheaf $\sA \subseteq \Sym^m \sB$ such that $\ka(\sA) \geq 1$.
  However, since $\sB$ is saturated in the locally free sheaf $\Omega_Y^1(\log
  D)$, it is reflexive, \cite[Claim on p.~158]{OSS}, and since
  $\Var(f^\circ)=1$, the sheaf $\sB$ is of rank $1$.  Thus $\sB \subseteq
  \Omega_Y^1(\log D)$ is an invertible subsheaf, \cite[Lem.~1.1.15, on
  p.~154]{OSS}, and Definition~\ref{def:speciallog} of a special pair gives
  that $\kappa(\sB) <1$, contradicting the fact that $\kappa(\sA) \geq 1$.  It
  follows that $\Var(f^\circ) = 0$ and that the family is hence isotrivial.
\end{proof}

A proof of Campana's Conjecture~\ref{conj:campana} in higher dimensions will
appear in a forthcoming paper.

\providecommand{\bysame}{\leavevmode\hbox to3em{\hrulefill}\thinspace}
\providecommand{\MR}{\relax\ifhmode\unskip\space\fi MR}
\providecommand{\MRhref}[2]{%
  \href{http://www.ams.org/mathscinet-getitem?mr=#1}{#2}
}
\providecommand{\href}[2]{#2}

\end{document}